\documentclass[a4paper]{article}
\usepackage{amsmath,amssymb,enumerate,bbm,theorem}

\newcommand{\TeXmacs}{T\kern-.1667em\lower.5ex\hbox{E}\kern-.125emX\kern-.1em\lower.5ex\hbox{\textsc{m\kern-.05ema\kern-.125emc\kern-.05ems}}}
\newcommand{\email}[1]{{\textit{Email:} \texttt{#1}}}
\newcommand{\mathd}{\mathrm{d}}
\newcommand{\tmem}[1]{{\em #1\/}}
\newcommand{\tmmathbf}[1]{\ensuremath{\boldsymbol{#1}}}
\newcommand{\tmop}[1]{\ensuremath{\operatorname{#1}}}
\newcommand{\tmstrong}[1]{\textbf{#1}}
\newcommand{\tmtextit}[1]{{\itshape{#1}}}
\newcommand{\withTeXmacstext}{This document has been produced using \TeXmacs (see \texttt{http://www.texmacs.org})}
\newenvironment{enumerateroman}{\begin{enumerate}[i.] }{\end{enumerate}}
\newenvironment{proof}{\noindent\textbf{Proof\ }}{\hspace*{\fill}$\Box$\medskip}

\theoremstyle{plain}
\newtheorem{theorem}{Theorem}[section]
\newtheorem{proposition}[theorem]{Proposition}
\newtheorem{corollary}[theorem]{Corollary}
\newtheorem{lemma}[theorem]{Lemma}

\theoremstyle{remark}
\newtheorem{remark}[theorem]{Remark}

\numberwithin{equation}{section}

\usepackage{hyperref}
\hypersetup{
pdfauthor = {Daniel GANDOLFO, Jean RUIZ, Marc WOUTS},
pdftitle={Limit Theorems and Coexistence Probabilities for the Curie-Weiss Potts Model with an external field},
pdfsubject = {},
pdfkeywords = {},
pdfcreator = {TeXmacs},
}

\begin{document}

\title{Limit Theorems and Coexistence Probabilities for the Curie-Weiss Potts
Model with an external field\thanks{{\withTeXmacstext}}}
\author{
Daniel GANDOLFO\thanks{\email{gandolfo@cpt.univ-mrs.fr}}\\
\small Centre de Physique Th\'eorique\\
\small UMR 6207 CNRS et Universit\'es Aix-Marseille et Sud Toulon-Var\\
\small Luminy Case 907, 13009 Marseille Cedex 9, France
\and
Jean RUIZ\thanks{\email{ruiz@cpt.univ-mrs.fr}}\\
\small Centre de Physique Th\'eorique\\
\small UMR 6207 CNRS et Universit\'es Aix-Marseille et Sud Toulon-Var\\
\small Luminy Case 907, 13009 Marseille Cedex 9, France
\and
Marc WOUTS\thanks{\email{wouts@math.univ-paris13.fr}}\\
\small Universit\'e Paris 13, CNRS, UMR 7539 LAGA\\
\small 99, avenue Jean-Baptiste Cl\'ement F-93 430 Villetaneuse France}
\maketitle

\begin{abstract}
  The Curie-Weiss Potts model is a mean field version of the well-known Potts
  model. In this model, the critical line $\beta = \beta_c (h)$ is explicitly
  known and corresponds to a first order transition when $q > 2$. In the
  present paper we describe the fluctuations of the density vector in the
  whole domain $\beta \geqslant 0$ and $h \geqslant 0$, including the
  conditional fluctuations on the critical line and the non-Gaussian
  fluctuations at the extremity of the critical line. The probabilities of
  each of the two thermodynamically stable states on the critical line are
  also computed. Similar results are inferred for the Random-Cluster model on
  the complete graph.
\end{abstract}

\section{Introduction}

The Curie-Weiss Potts model is a model of statistical mechanics which, being a
mean-field model, can be studied by means of analytic tools. First it was
shown in {\cite{Wu}} that at $h = 0$, the model undergoes a phase transition
at the critical inverse temperature
\[ \beta_c = \left\{ \begin{array}{ll}
     q & \text{if } q \leqslant 2\\
     2 \frac{q - 1}{q - 2} \log (q - 1) & \text{if } q > 2.
   \end{array} \right. \]
When $q > 2$ this transition is first order. The case of non-zero external
field was considered in {\cite{BCC}} and it appeared that the first-order
transition remains on a critical line. Recently this critical line was
computed explicitly {\cite{BGRW}}.

On the critical line, two or more states can coexist. One of the issue we
address in the present work is the computation of the probabilities of these
stable states. We also obtain a description of the limit distribution of the
empirical vector of the spin variables that extend previous results on the
Curie-Weiss Ising model {\cite{EN-JSP}} (see also {\cite{EN-PTRF,ENR}}), and
previous results on the Curie-Weiss Potts model with no external field
{\cite{EW}}.

The Curie-Weiss Potts model is connected as well to the random-cluster model.
In that model, the first order phase transition for $q > 2$ was described in
{\cite{BGJ}} and it appeared that at criticality, two possible structures of
the random graph are possible. The probability for each structure was latter
computed in {\cite{LL}}. A consequence of our results we present a simple way
of computing these probabilities when $q > 2$ is integer.

\section{The Curie-Weiss Potts model}

The Curie-Weiss Potts model is a spin model on the complete graph. The
probability of observing the configuration $\sigma \in \{1, \ldots, q\}^n$ at
inverse temperature $\beta$, in an exterior field $H = h / \beta$ equals
\[ \mu_{\beta, h, n} (\sigma) = \frac{1}{Z_{\beta, h, n}} \exp \left(
   \frac{\beta}{n}  \sum_{1 \leqslant i < j \leqslant n} \delta_{\sigma_i,
   \sigma_j} + h \sum_{i = 1}^n \delta_{\sigma_i, 1} \right) \]
where $\delta$ is the Kronecker symbol and $Z_{\beta, h, n}$ the partition
function
\[ Z_{\beta, h, n} = \sum_{\sigma \in \{1, \cdots, q\}^n} \exp \left(
   \frac{\beta}{n}  \sum_{1 \leqslant i < j \leqslant n} \delta_{\sigma_i,
   \sigma_j} + h \sum_{i = 1}^n \delta_{\sigma_i, 1} \right) . \]
Our interest is in the limit distribution of the empirical vector
\begin{equation}
  \tmmathbf{N}= (N_1, \ldots, N_q) = \left( \sum_{i = 1}^n \delta_{\sigma_i,
  1}, \ldots, \sum_{i = 1}^n \delta_{\sigma_i, q} \right)
\end{equation}
that represents the number of spins of each color for a given configuration
$\sigma$. The normalized vector $\tmmathbf{N}/ n$ belongs to the set of
probability vectors
\begin{equation}
  \Omega =\{\tmmathbf{x} \in \mathbbm{R}^q : x_1 + \cdots + x_q = 1 \text{ and
  } x_i \geqslant 0, \forall i\}. \label{defOm}
\end{equation}
The large deviation principle for $\tmmathbf{N}/ n$ is an immediate
application of Stirling's formula (see for instance Lemma \ref{lem-stirling}).
If we consider $f_{\beta, h}$ the microcanonical free energy of the model:
\begin{equation}
  f_{\beta, h} (\tmmathbf{x}) = \sum_{i = 1}^q x_i \log x_i - \frac{\beta}{2}
  \sum_{i = 1}^q x_i^2 - hx_1 \text{, \ \ } \forall \tmmathbf{x} \in \Omega
  \label{deff}
\end{equation}
with the convention that $0 \log 0 = 0$, then we have the following classical
large deviation result (see for instance {\cite{CET}}, and also
{\cite{BGJ,BCS}} for LDP concerning the closely related random cluster model).

\begin{theorem}
  \label{thm-LDP}Assume that $\beta_n \rightarrow \beta$ and $h_n \rightarrow
  h$. Then, the vector $\tmmathbf{N}/ n \in \Omega$ distributed according to
  the measure $\mu_{\beta_n, h_n, n}$ follows a large deviation principle with
  speed $n$ and good rate function $f_{\beta, h} - \min_{\Omega} f_{\beta,
  h}$.
\end{theorem}

This large deviation principle leads to a law of large number: when $f_{\beta,
h}$ has a unique global minimizer, $\tmmathbf{N}/ n$ converges towards that
minimizer. The structure of the minimizers of $f_{\beta, h}$ was determined in
the papers {\cite{Wu,KS,BCC}}. Here we give some further details:

\begin{proposition}
  \label{prop-struct-x}Let $\beta, h \geqslant 0$ and let $\tmmathbf{x}$ be a
  global minimizer of $f_{\beta, h}$ in $\Omega$.
  \begin{enumerateroman}
    \item The vector $\tmmathbf{x}$ has the coordinate $\min (x_i)$ repeated
    $q - 1$ times at least.
    
    \item If $h > 0$, then $x_1 > x_i$, for all $i \in \{2, \ldots, q\}$.
    
    \item The inequality $\min (x_i) > 0$ holds.
    
    \item For any $q \geqslant 3$, or $q = 2$ and $(\beta, h) \neq (\beta_c,
    0)$, one has $\min (x_i) < 1 / \beta$.
  \end{enumerateroman}
\end{proposition}

Because of the simple structure of the global minimizers of the free energy,
the problem of finding them reduces to a one-dimensional optimization problem.
The usual parametrization consists in taking $x_1 = (1 + (q - 1) s) / q$, $x_2
= \cdots = x_q = (1 - s) / q$ where $s \in [0, 1]$ is a parameter called the
magnetization. Another equivalent parametrization permitted in {\cite{BGRW}}
the explicit computation of the critical line

\begin{equation}
  h_T = \left\{ (\beta, h) : 0 \leqslant h < h_0 \text{ and } h = \log (q - 1)
  - \beta \frac{q - 2}{2 (q - 1)} \right\}
\end{equation}
with extremities $(\beta_c, 0)$ and $(\beta_0, h_0)$, where
\[ \beta_0 = 4 \frac{q - 1}{q} \text{ \ and \ } h_0 = \log (q - 1) - 2 \frac{q
   - 2}{q} \]
were already determined in {\cite{BCC}}. The key observation in {\cite{BGRW}}
was that the free energy $f_{\beta, h} (\tmmathbf{x}_z)$ at
\begin{equation}
  \tmmathbf{x}_z = \left( \frac{1 + z}{2}, \frac{1 - z}{2 (q - 1)}, \cdots,
  \frac{1 - z}{2 (q - 1)} \right) \label{xz} \text{, \ \ \ } z \in [\pm 1]
\end{equation}
is easily split into its even and odd parts:
\begin{eqnarray*}
  f_{\beta, h} (\tmmathbf{x}_z) & = & \frac{1 + z}{2} \log \frac{1 + z}{2} +
  \frac{1 - z}{2} \log \frac{1 - z}{2} - \frac{1}{2} \log (q - 1) -
  \frac{\beta (1 + z^2)}{8} \left[ 1 + \frac{1}{q - 1} \right] - \frac{1}{2}
  h\\
  &  & + \frac{z}{2} \left[ \log (q - 1) - \beta \frac{q - 2}{2 (q - 1)} - h
  \right]
\end{eqnarray*}
showing that, on the critical line $h_T$, the free energy $f_{\beta, h}
(\tmmathbf{x}_z)$ is an {\tmem{even}} function of $z$. It is strictly convex
for $\beta < \beta_0$ but not for $\beta \geqslant \beta_0$. Indeed, the
second derivative of $z \mapsto f_{\beta, h} (\tmmathbf{x}_z)$ is
\begin{eqnarray}
  \frac{\mathd^2 f_{\beta, h} (\tmmathbf{x}_z)}{\mathd z^2} & = & \frac{1}{1 -
  z^2} - \frac{\beta q}{4 (q - 1)},  \label{d2fz}
\end{eqnarray}
thus, for $\beta \geqslant \beta_0$, the function is strictly convex on $[- 1,
- z_i)$ and on $(z_i, 1]$, concave on $(- z_i, z_i)$ where
\begin{equation}
  z_i = \sqrt{1 - \beta_0 / \beta} .
\end{equation}
Depending on the parameters $(\beta, h)$ the free energy presents one or
several global minimizers. The following is a summary of the works {\cite{Wu}}
(for $h = 0$) and {\cite{BGRW}} (for $h > 0$):

\begin{theorem}
  Let $\beta, h \geqslant 0$.\label{thm-transition}
  \begin{enumerateroman}
    \item If $h > 0$ and $(\beta, h) \notin h_T$, the free energy $f_{\beta,
    h}$ has a unique global minimizer in $\Omega$. This minimizer is analytic
    in $\beta$ and $h$ outside of $h_T \cup \{(\beta_0, h_0)\}$.
    
    \item If $h > 0$ and $(\beta, h) \in h_T$, the free energy $f_{\beta, h}$
    has two global minimizers in $\Omega$. More precisely, for any $z \in (0,
    (q - 2) / q)$, the two global minimizers of $f_{\beta_z, h_z}$ at
    \[ \beta_z = 2 \frac{q - 1}{q}  \frac{1}{z} \log \frac{1 + z}{1 - z}
       \text{ \ \ and \ \ } h_z = \log (q - 1) - \frac{q - 2}{2 (q - 1)}
       \beta_z \]
    are the points $\tmmathbf{x}_{\pm z}$. Furthermore, $\tmmathbf{x}_z$
    (resp. $\tmmathbf{x}_{- z}$) is the limit of the unique global minimizer
    of $f_{\beta, h}$ as $(\beta, h) \rightarrow (\beta_z, h_z)$ above (resp.
    below) the line $h_T$.
    
    \item If $h = 0$ and $\beta < \beta_c$, the unique global minimizer of
    $f_{\beta, h}$ is $(1 / q, \ldots, 1 / q) =\tmmathbf{x}_{- (q - 2) / q}$.
    
    \item If $h = 0$ and $\beta > \beta_c$, there are $q$ global minimizers of
    $f_{\beta, h}$, which all equal $\tmmathbf{x}_z$ up to a permutation of
    the coordinates, for some appropriate $z \in ((q - 2) / q, 1)$.
    
    \item If $h = 0$ and $\beta = \beta_c$, there are $q + 1$ global
    minimizers of $f_{\beta, h}$ : the symmetric one $(1 / q, \ldots, 1 / q)
    =\tmmathbf{x}_{- (q - 2) / q}$ together with the permutations of
    \[ \left( \frac{q - 1}{q}, \frac{1}{q (q - 1)}, \cdots, \frac{1}{q (q -
       1)} \right) =\tmmathbf{x}_{(q - 2) / q} . \]
  \end{enumerateroman}
\end{theorem}

\section{Statement of the results}

In this paper we address essentially two questions. According to Theorem
\ref{thm-LDP} the distribution of $\tmmathbf{N}/ n$ is concentrated, as $n
\rightarrow + \infty$, on the set of global minimizers of the free energy.
First, we study the fluctuations of the empirical vector $\tmmathbf{N}$ around
its typical value. Second, when several global minimizers exists we explicit
the weight of each of them.

These questions were answered in several very interesting papers for
particular cases of the model. The case of the Curie-Weiss Ising model ($q =
2$) was reported in {\cite{EN-JSP}} (see {\cite{EN-PTRF,ENR}} for the proofs),
while the Curie-Weiss Potts model was treated at zero external field in
{\cite{EW}}.

Our approach is similar to that of the former references, with the technical
difference that our computations are based on Stirling's formula while the
former works are based on the fact that the law of $\tmmathbf{N}/ n
+\tmmathbf{W}/ \sqrt{n}$, where $\tmmathbf{W}$ is a Gaussian vector in
$\mathbbm{R}^q$ with distribution $\mathcal{N}(0, \beta^{- 1} I_q)$, can be
explicitly computed (see for instance Lemma 3.2 in {\cite{EW}}).

We also permit that the parameters $\beta$ and $h$ fluctuate with $n$, and
take in the sequel $(\beta_n, h_n) \rightarrow (\beta, h)$. This will be
useful for applying our results to related model such as the random cluster
model on the complete graph.

Our first result concerns the fluctuations of the empirical vector
$\tmmathbf{N}$ outside of the critical line. The fluctuations belong to the
hyperplane
\begin{eqnarray}
  \mathcal{H} & = & \left\{ \tmmathbf{w} \in \mathbbm{R}^d : \sum_{i = 1}^d
  w_i = 0 \right\} .  \label{def-hypH}
\end{eqnarray}
Not surprisingly, these fluctuations are Gaussian. This generalizes Theorem
2.4 in {\cite{EW}} to the case of positive external fields. The way that
$(\beta_n, h_n)$ converges to $(\beta, h)$ is able to shift the center of the
distribution.

\begin{theorem}
  \label{thm-fluct1}Assume that $(\beta_n, h_n) \rightarrow (\beta, h)$ for
  some $\beta, h \geqslant 0$ with $(\beta, h) \neq (\beta_0, h_0)$. Assume
  that there is a unique global minimizer $\tmmathbf{x}= (x_1, x_q, \ldots,
  x_q)$ of the free energy $f_{\beta, h}$. For every $n$, let $\tmmathbf{d}_n$
  the smallest $\tmmathbf{d} \in \mathcal{H}$ such that
  $\tmmathbf{x}+\tmmathbf{d} \in \Omega$ is a local minimizer of $f_{\beta_n,
  h_n}$. Let $\tmmathbf{W}$ be the random variable in $\mathcal{H}$ such that
  \begin{equation}
    \tmmathbf{N}= n\tmmathbf{x}+ n\tmmathbf{d}_n + n^{1 / 2} \tmmathbf{W}
    \label{defW}
  \end{equation}
  where the distribution of $\tmmathbf{N}$ is given according to the measure
  $\mu_{\beta_n, h_n, n}$. Then, {\tmstrong{W}} converges in law towards the
  centered Gaussian vector with covariance matrix
  \begin{equation}
    \left( \frac{1}{x_1 x_q} - q \beta \right)^{- 1} \left(\begin{array}{cccc}
      q - 1 & - 1 & \cdots & - 1\\
      - 1 & 1 + (q - 2) \frac{\frac{1}{x_1} - \beta}{\frac{1}{x_q} - \beta} & 
      & - \frac{\frac{1}{x_1} - \beta}{\frac{1}{x_q} - \beta}\\
      \vdots &  & \ddots & \\
      - 1 & - \frac{\frac{1}{x_1} - \beta}{\frac{1}{x_q} - \beta} &  & 1 + (q
      - 2) \frac{\frac{1}{x_1} - \beta}{\frac{1}{x_q} - \beta}
    \end{array}\right) \label{defKx}
  \end{equation}
  which has rank $q - 1$.
\end{theorem}

\begin{remark}
  \label{rmk-d}The vector $\tmmathbf{d}_n$ is $O \left( | \beta_n - \beta | +
  |h_n - h| \right)$ when the quadratic term in the Taylor expansion of
  $f_{\beta, h}$ is definite, that is for any $(\beta, h) \neq (\beta_0, h_0)$
  -- see Lemmas \ref{lem-Taylor-x}, \ref{lem-Taylor-bh} and
  \ref{lem-quad-form} below. Hence, for $\beta_n - \beta = o (n^{- 1 / 2})$
  and $h_n - h = o (n^{- 1 / 2})$ the vector $n\tmmathbf{d}_n$ is negligible
  with respect to $n^{1 / 2} \tmmathbf{W}$ and could be removed from the
  definition of $\tmmathbf{W}$ at (\ref{defW}). It is remarkable also that on
  the line $\beta < \beta_c$, $h = 0$ we have $\tmmathbf{x}= (1 / q, \ldots, 1
  / q)$, hence for $h_n = 0$ and $\beta_n \rightarrow \beta$, the vector
  $\tmmathbf{d}_n$ is exactly zero.
\end{remark}

\begin{remark}
  In the range of validity of Theorem 2.4 in {\cite{EW}}, that is $\beta_n =
  \beta < \beta_c$ and $h_n = h = 0$, we have $x_1 = \ldots = x_q = 1 / q$
  thus the covariance matrix simplifies to
  \[ \frac{1}{q^2 - q \beta}  \left(\begin{array}{ccc}
       q - 1 &  & - 1\\
       & \ddots & \\
       - 1 &  & q - 1
     \end{array}\right) . \]
  We have checked the correspondence with the covariance matrix that appears
  in {\cite{EW}}.
\end{remark}

The matrix (\ref{defKx}) gives a special emphasis on the first coordinate
since it corresponds to the case $x_2 = \ldots = x_q$. Before stating the next
theorem we give a more symmetric definition for the covariance matrix: we let
\begin{eqnarray}
  K (\tmmathbf{x}) & = & \left( \frac{1}{\min (x_i) \max (x_i)} - q \beta
  \right)^{- 1}  \left(\begin{array}{ccc}
    1 + (q - 2) \alpha (x_1, x_1) &  & - \alpha (x_1, x_q)\\
    & \ddots & \\
    - \alpha (x_q, x_1) &  & 1 + (q - 2) \alpha (x_q, x_q)
  \end{array}\right)  \label{def-Kx}
\end{eqnarray}
where
\[ \alpha (x, y) = \frac{\max (x_i)^{- 1} - \beta}{\max (x, y)^{- 1} - \beta}
   . \]

When the free energy has several global minimizers, that is when $(\beta, h)
\in h_T$ or $\beta \geqslant \beta_c$ and $h = 0$, the empirical vector
$\tmmathbf{N}/ n$ is close to either one or the other of the minimizers of the
free energy $f_{\beta, h}$. We first determine the conditional fluctuations
(this extends Theorem 2.5 of {\cite{EW}}):

\begin{theorem}
  \label{thm-fluct-cond}Assume that $(\beta_n, h_n) \rightarrow (\beta, h)$
  with $\beta, h \geqslant 0$. Assume that the free energy $f_{\beta, h}$ has
  multiple global minimizers $\tmmathbf{x}, \tmmathbf{x}', \ldots$ and let
  $\varepsilon > 0$ smaller than the distance between any two global
  minimizers of $f_{\beta, h}$. Let $\tmmathbf{d}_n$ the smallest
  $\tmmathbf{d} \in \mathcal{H}$ such that $\tmmathbf{x}+\tmmathbf{d} \in
  \Omega$ is a local minimizer of $f_{\beta_n, h_n}$. Then, under the
  conditional measure
  \[ \mu_{\beta_n, h_n, n} \left( . \left| \frac{\tmmathbf{N}}{n} \in B
     (\tmmathbf{x}, \varepsilon) \right. \right), \]
  the variable $\tmmathbf{W}$ defined by $\tmmathbf{N}= n\tmmathbf{x}+
  n\tmmathbf{d}_n + n^{1 / 2} \tmmathbf{W}$ converges in law to the centered
  Gaussian vector with covariance matrix $K (\tmmathbf{x})$.
\end{theorem}

Additionally we compute the limit probabilities that $\tmmathbf{N}/ n$ be
close to a given global minimizer of the free energy, generalizing Theorem 2.3
of {\cite{EW}} with an explicit formula.

\begin{theorem}
  \label{thm-probacoex}Assume that there are $\beta, h \geqslant 0$ and
  $\lambda, \nu \in \mathbbm{R}$ such that
  \[ (\beta_n, h_n) = (\beta, h) + n^{- 1} \left( \lambda, \nu \right) + o
     (n^{- 1}), \]
  and assume that the free energy $f_{\beta, h}$ has multiple global
  minimizers $\tmmathbf{x}, \tmmathbf{x}', \ldots$ If $\varepsilon > 0$ is
  smaller than the distance between any two global minimizers of $f_{\beta,
  h}$, then
  \begin{eqnarray}
    \lim_{n \rightarrow \infty} \mu_{\beta_n, h_n, n} \left(
    \frac{\tmmathbf{N}}{n} \in B (\tmmathbf{x}, \varepsilon) \right) & = &
    \frac{\tau (\tmmathbf{x})}{\tau (\tmmathbf{x}) + \tau (\tmmathbf{x}') +
    \cdots}  \label{eq-probas}
  \end{eqnarray}
  where
  \begin{eqnarray}
    \tau (\tmmathbf{x}) & = & \left( 1 - \beta \min_{i = 1}^q (x_i)
    \right)^{\frac{2 - q}{2}} \exp \left( \frac{\lambda}{2} \sum_{i = 1}^q
    x_i^2 + \nu x_1 \right) . 
  \end{eqnarray}
\end{theorem}

\begin{remark}
  On the critical line $h_T$ one can parametrize the formula (\ref{eq-probas})
  according to the second point in Theorem \ref{thm-transition} : when
  $(\beta, h) = (\beta_z, h_z)$ with $z \in (0, (q - 2) / q)$ the two global
  minimizers are $\tmmathbf{x}_{\pm z}$. In particular, when $z \rightarrow 0$
  (i.e. $(\beta, h)$ on $h_T$ close to $(\beta_0, h_0)$), the probability of
  each corresponding state converges to $1 / 2$.
\end{remark}

We also describe the fluctuations at the extremity $(\beta_0, h_0)$ of the
critical line. This extends for instance Theorem 2 in {\cite{EN-JSP}} that
applies to the case of the Curie-Weiss Ising model at criticality, namely $q =
2$ and $(\beta_0, h_0) = (\beta_c, 0)$. We recall that $\mathcal{H}$ defined
at (\ref{def-hypH}) is the hyperplane parallel to $\Omega$. Given a vector
$\tmmathbf{u} \in \mathbbm{R}^q$, we denote by $\tmmathbf{u}^{\perp}$ the
vector space made of all vectors orthogonal to $\tmmathbf{u}$ in the Euclidean
space $\mathbbm{R}^q$.

\begin{theorem}
  \label{thm-fluct-b0h0}Assume that $(\beta_n, h_n) \rightarrow (\beta_0,
  h_0)$ with $\beta_n - \beta = o (n^{- 3 / 4})$ and $h_n - h = o (n^{- 3 /
  4})$ and let $\tmmathbf{x}= (1 / 2, 1 / 2 (q - 1), \ldots, 1 / 2 (q - 1))$
  be the unique minimizer of $f_{\beta_0, h_0}$. Let $\tmmathbf{u}= (1 - q, 1,
  \ldots, 1)$. If the random variables $T \in \mathbbm{R}$ and $\tmmathbf{V}
  \in \mathcal{H} \cap \tmmathbf{u}^{\perp}$ are defined by
  \begin{equation}
    \tmmathbf{N}= n\tmmathbf{x}+ n^{3 / 4} T\tmmathbf{u}+ n^{1 / 2}
    \tmmathbf{V} \label{defTV},
  \end{equation}
  then $(T, \tmmathbf{V})$ converges in law. The limit has the following
  properties:
  \begin{enumerateroman}
    \item $T$ and $\tmmathbf{V}$ are asymptotically independent
    
    \item $T$ converges in law to the probability measure on $\mathbbm{R}$
    proportional to
    \[ \exp \left( - \frac{4 (q - 1)^4}{3} t^4 \right) \mathd t \]
    \item $\tmmathbf{V}$ converges in law towards the centered Gaussian vector
    with covariance matrix
    \[ \frac{q}{2 (q - 1)^2 (q - 2)}  \left(\begin{array}{cccc}
         0 & 0 & \cdots & 0\\
         0 & q - 2 &  & - 1\\
         \vdots &  & \ddots & \\
         0 & - 1 &  & q - 2
       \end{array}\right) \]
    of rank $q - 2$.
  \end{enumerateroman}
\end{theorem}

We conclude the summary of our results with two claims on the random-cluster
model $G (n, p, q)$ on the complete graph $K_n$ with $n$ vertices. In that
model, a configuration $\omega \in \{0, 1\}^{E (K_n)}$ has a probability
proportional to
\[ \prod_{e \in E (K_n)} p^{\omega_e} (1 - p)^{1 - \omega_e} q^{C (\omega)} \]
where $C (\omega)$ stands for the number of connected components of the
sub-graph with edge set $\{e \in E (K_n) : \omega_e = 1\}$. This model is
closely related to the Potts model after the well known Fortuin-Kasteleyn
representation (see for instance {\cite{ES}}). We take a spin configuration
$\sigma \in \{1, \ldots, q\}^{V (K_n)}$ under the measure $\mu_{\beta_n, h_n,
n}$, then let $\omega_e = 1$ with probability $p_n = 1 - \exp (- \beta_n / n)$
only if $\sigma_i = \sigma_j$, where $i, j$ are the extremities of the edge
$e$ (and else $\omega_e = 0$). The resulting configuration $\omega$ follows
the distribution of the random cluster model $G (n, p_n, q)$.

First we have a Corollary of Theorem \ref{thm-probacoex}: we compute the
probability that there exists a giant component in $G (n, p_n, q)$, that is a
connected component for $\omega$ of size $\Theta (n)$, when $p_n$ is close to
the critical value $\beta_c / n$. This completes part (b) of Theorem 2.3 in
{\cite{BGJ}}, with a simpler proof than that of Theorem 19 of {\cite{LL}}.

\begin{corollary}
  \label{cor-proba-gc}Let $q > 2$ integer and consider $p_n$ such that
  \[ p_n = \frac{\beta_c}{n} + \frac{\gamma}{n^2} + o \left( \frac{1}{n^2}
     \right) . \]
  Then, with a probability that converges to
  \[ \frac{1}{1 + \frac{1}{q}  \left( \frac{1 - \beta_c / q}{1 - \beta_c / (q
     (q - 1))} \right)^{\frac{2 - q}{2}} \exp \left( - \left(
     \frac{\beta_c^2}{4} + \frac{\gamma}{2} \right)  \frac{(q - 2)^2}{q (q -
     1)} \right)} \]
  the graph $G (n, p_n, q)$ contains a giant component.
\end{corollary}

The description of the Gaussian fluctuations also enable fine computations of
the partition function of the random-cluster model
\begin{eqnarray}
  Z^{\tmop{RC}}_{p, q, n} & = & \sum_{\omega \in \{0, 1\}^{E (K_n)}} \prod_{e
  \in E (K_n)} p^{\omega_e} (1 - p)^{1 - \omega_e} q^{C (\omega)} . 
  \label{eq-def-ZRC}
\end{eqnarray}
For instance,

\begin{proposition}
  \label{prop-Zrc}The partition function of the random cluster model for
  integer $q \geqslant 2$ and
  \[ p_n = \frac{\beta}{n} + \frac{\gamma}{n^2} + o \left( \frac{1}{n^2}
     \right) \]
  with $0 \leqslant \beta < \beta_c$ and $\gamma \in \mathbbm{R}$ satisfies
  \begin{eqnarray}
    Z^{\tmop{RC}}_{p_n, q, n} & = & (1 + o_n (1)) \left( 1 - \frac{\beta}{q}
    \right)^{- \frac{q - 1}{2}} q^n \exp \left( - \frac{2 n \beta + 2 \gamma +
    \beta^2}{4} \left( \frac{q - 1}{q} \right) \right) .  \label{eq-ZRC}
  \end{eqnarray}
\end{proposition}

\begin{remark}
  Although our Theorem \ref{thm-probacoex} agrees with Theorem 2.3 of
  {\cite{EW}} when $h = 0$, Corollary \ref{cor-proba-gc} and Proposition
  \ref{prop-Zrc} do not give exactly the same conclusions as, respectively,
  Theorem 19 and Theorem 9 (i) in {\cite{LL}}. The latter Theorem states an
  equivalent to the partition function restricted to the set of configurations
  made of trees and unicyclic components, which, for $\beta < \beta_c$, is
  equivalent to the whole partition function. The ratio of the equivalent in
  Theorem 9 (i) in {\cite{LL}} over (\ref{eq-ZRC}) is
  \[ \exp \left( - \frac{3}{4} + \frac{\beta}{2} + \frac{\beta^2}{4 q} \right)
  \]
  (the formulas do coincide at the exponential order). We could not check the
  proofs in {\cite{LL}}, yet we were surprised to find that Theorem 9 (i)
  would not permit to recover $Z^{\tmop{RC}}_{p_n, 1, n} = 1$ for $q = 1$.
\end{remark}

\section{Proofs}

\label{sec-proofs}This section is organized as follows. First we describe the
asymptotics of the distribution using Stirling's formula. We also prove
Proposition \ref{prop-struct-x}. Then we address successively the limit
distribution at $(\beta, h) \neq (\beta_0, h_0)$, at the extremity $(\beta_0,
h_0)$ of the critical line, and finally we give the proofs of the related
results for the random-cluster model.

\subsection{Asymptotic density \& limit of the uniform measure}

In this Section we give an equivalent to the density of the Potts model, prove
Proposition \ref{prop-struct-x} and describe the limit of the uniform measure
on the set of possible realizations of $\tmmathbf{N}/ n$.

For any $\varepsilon \geqslant 0$, we let
\[ \Omega_{\varepsilon} =\{\tmmathbf{x} \in \Omega : \min x_i \geqslant
   \varepsilon\} \text{ \ \ and \ \ } \Omega_{\varepsilon}^n =\{\tmmathbf{x}
   \in \Omega_{\varepsilon} : n\tmmathbf{x} \in \mathbbm{N}^q \}, \]
where $\Omega$ is the set of probability vectors, see (\ref{defOm}). We also
write $\Omega_{0^+} = \bigcup_{\varepsilon > 0} \Omega_{\varepsilon}$ and
$\Omega^n = \Omega_0^n$. In our first Lemma we give an equivalent to the
density of the Potts model with respect to the counting measure on $\Omega^n$.
We use nothing else than Stirling's formula
\begin{eqnarray*}
  n! & = & (1 + o_n (1)) \sqrt{2 \pi n}  \left( \frac{n}{e} \right)^n .
\end{eqnarray*}
For any $\tmmathbf{x} \in \Omega_{0^+}$ and $\beta \geqslant 0$, we let
\begin{equation}
  \left. A_{\beta} (\tmmathbf{x}) = (2 \pi \right)^{- \frac{q - 1}{2}} 
  \prod_{i = 1}^q x_i^{- 1 / 2} \exp \left( - \frac{\beta}{2} \right) .
\end{equation}
We also recall that the free energy $f_{\beta, h}$ was defined at
(\ref{deff}). We have:

\begin{lemma}
  \label{lem-stirling}For any $\beta, h$, any $n \geqslant 1$ and
  $\tmmathbf{x} \in \Omega_{0^+}^n$, define $r_{\beta, h, n} (\tmmathbf{x})$
  by
  \[ Z_{\beta, h, n} \mu_{\beta, h, n} \left( \tmmathbf{N}= n\tmmathbf{x}
     \right) = (1 + r_{\beta, h, n} (\tmmathbf{x})) n^{- \frac{q - 1}{2}}
     A_{\beta} (\tmmathbf{x}) \exp (- nf_{\beta, h} (\tmmathbf{x})) . \]
  Then, for any $\varepsilon > 0$, $\sup_{x \in \Omega_{\varepsilon}^n}
  \sup_{\beta, h} |r_{\beta, h, n} (\tmmathbf{x}) |$ goes to $0$ as $n
  \rightarrow \infty$.
\end{lemma}

\begin{proof}
  Given $\tmmathbf{x} \in \Omega_{\varepsilon}^n$ we write $n\tmmathbf{x}=
  (n_1, \ldots, n_q) =\tmmathbf{n}$. It is a vector with positive integer
  coordinates. There are exactly $n! / \prod_{i = 1}^q n_i$ ways of choosing
  the spin configuration that satisfy the constraint $\left( N_1, \ldots, N_q
  \right) = (n_1, \ldots, n_q)$, hence
  \begin{eqnarray*}
    Z_{\beta, h, n} \mu_{\beta, h, n} \left( \tmmathbf{N}=\tmmathbf{n} \right)
    & = & \frac{n!}{n_1 ! \cdots n_q !} \exp \left( \frac{\beta}{n} \sum_{i =
    1}^q \frac{n_i (n_i - 1)}{2} + hn_1 \right)\\
    & = & \frac{n!}{n_1 ! \cdots n_q !} \exp \left( - \frac{\beta}{2} + n
    \left[ \frac{\beta}{2} \sum_{i = 1}^q x_i^2 + hx_1 \right] \right) .
  \end{eqnarray*}
  Thus
  \begin{eqnarray*}
    1 + r_{\beta, h, n} (\tmmathbf{x}) & = & \frac{Z_{\beta, h, n} \mu_{\beta,
    h, n} \left( \tmmathbf{N}= n\tmmathbf{x} \right)}{n^{- \frac{q - 1}{2}}
    A_{\beta} (\tmmathbf{x}) \exp (- nf_{\beta, h} (\tmmathbf{x}))}\\
    & = & \frac{n!}{n_1 ! \cdots n_q !}  \frac{\prod_{i = 1}^q \sqrt{2 \pi
    nx_i}}{\sqrt{2 \pi n}} \exp \left( n \sum_{i = 1}^q x_i \log x_i \right)
  \end{eqnarray*}
  which does not depend on $\beta$ nor on $h$. Applying Stirling's formula
  yields the conclusion as all the $n_i$ go to infinity uniformly over
  $\tmmathbf{x} \in \Omega_{\varepsilon}^n$.
\end{proof}

\begin{remark}
  Theorem \ref{thm-LDP} is a consequence of Lemma \ref{lem-stirling} as
  $f_{\beta_n, h_n} (\tmmathbf{x}) \underset{n \rightarrow
  \infty}{\longrightarrow} f_{\beta, h} (\tmmathbf{x})$ when $(\beta_n, h_n)
  \rightarrow (\beta, h)$, uniformly over $\tmmathbf{x} \in \Omega$ (for the
  uniformity, see Lemma \ref{lem-Taylor-bh} below).
\end{remark}

Now we give the proof of Proposition \ref{prop-struct-x}. For completeness we
repeat some arguments from {\cite{KS,BCC}}.

\begin{proof}
  (Proposition \ref{prop-struct-x}). Call $g_{\beta} (z) = \log (z) - \beta
  z$. As $\tmmathbf{x}$ is a minimizer of the free energy in $\mathcal{H}$ one
  has
  \begin{eqnarray}
    0 & = & \frac{\partial f_{\beta, h}}{\partial x_i} (\tmmathbf{x}) -
    \frac{\partial f_{\beta, h}}{\partial x_j} (\tmmathbf{x}) \nonumber\\
    & = & g_{\beta} (x_i) - g_{\beta} (x_j) - h \left[ \tmmathbf{1}_{\{i =
    1\}} -\tmmathbf{1}_{\{j = 1\}} \right] \text{, \ \ \ \ } \forall i, j \in
    \{1, \ldots q\}  \label{eq-g}\\
    \text{and \ } 0 & \leqslant & \frac{\partial^2 f_{\beta, h}}{\partial
    x_i^2} (\tmmathbf{x}) + \frac{\partial^2 f_{\beta, h}}{\partial x_j^2}
    (\tmmathbf{x}) \nonumber\\
    & = & g'_{\beta} (x_i) + g'_{\beta} (x_j) .  \label{eq-dg}
  \end{eqnarray}
  First we assume $h = 0$. As $g$ is concave, (\ref{eq-g}) implies that the
  set $\{x_i : i = 1, \ldots, q\}$ contains at most two values. Equation
  (\ref{eq-dg}) implies that at most one of the $x_i$ has $g'_{\beta} (x_j) <
  0$. As
  \begin{eqnarray*}
    g'_{\beta} (x) & = & \frac{1}{x} - \beta
  \end{eqnarray*}
  is positive on $(0, 1 / \beta)$ and negative on $(1 / \beta, 1)$, the first
  point of Proposition \ref{prop-struct-x} follows together with the
  inequality
  \begin{eqnarray}
    \min (x_i) & \leqslant & \frac{1}{\beta} .  \label{eq-minxi-l}
  \end{eqnarray}
  Assume now that $h > 0$ and that (\ref{eq-minxi-l}) does not hold for some
  $i \in \{2, \ldots, q\}$. If $x_i > x_1$ with $i \in \{2, \ldots, q\}$, the
  vector $\tilde{\tmmathbf{x}}$ with $x_1$ and $x_i$ permuted has $f_{\beta,
  h} ( \tilde{\tmmathbf{x}}) < f_{\beta, h} (\tmmathbf{x})$, a contradiction,
  therefore $x_1 \geqslant x_i$. The equality $x_1 = x_2$ is impossible in
  view of (\ref{eq-g}), yielding the second point of Proposition
  \ref{prop-struct-x}. Now we conclude the proof of the first point of
  Proposition \ref{prop-struct-x} when $h > 0$ : the inequality $x_i \geqslant
  1 / \beta$, that implies $g'_{\beta} (x_i) \leqslant 0$ and $g'_{\beta}
  (x_1) < 0$ since $x_1 > x_i$, would contradict (\ref{eq-dg}). Hence all the
  $x_i$ belong to $(0, 1 / \beta)$ where there is at most one reciprocal image
  of $g_{\beta} (x_1) - h$ by $g_{\beta}$, hence $x_2 = \ldots = x_q < 1 /
  \beta$.
  
  Now we address the third point. If $\min (x_i) = 0$, one can find $i, j \in
  \{1, \ldots, q\}$ such that $x_i = 0$ and $x_j > 0$ as $\tmmathbf{x} \in
  \Omega$. Hence $\tmmathbf{x}^t =\tmmathbf{x}+ t (\tmmathbf{e}_i
  -\tmmathbf{e}_j)$ belongs to $\Omega_{0^+}$ for small enough $t > 0$. Yet,
  \begin{eqnarray*}
    \frac{\mathd}{\mathd t} f_{\beta, h} (\tmmathbf{x}^t) & = & g_{\beta} (x_i
    + t) - g_{\beta} (x_j - t) + h \left[ \tmmathbf{1}_{\{i = 1\}}
    -\tmmathbf{1}_{\{j = 1\}} \right]
  \end{eqnarray*}
  goes to $- \infty$ as $t \rightarrow 0^+$, a contradiction.
  
  Remains the strict inequality in (\ref{eq-minxi-l}). We let
  $\tmmathbf{x}^{\beta} = (1 - (q - 1) / \beta, 1 / \beta, \ldots, 1 / \beta)$
  (which is in $\Omega$ for $\beta \geqslant q - 1$, and satisfies the case of
  equality in (\ref{eq-minxi-l}) for $\beta \geqslant q$) and derive
  conditions for $\tmmathbf{x}_{\beta}$ being a minimizer of the free energy.
  Equation (\ref{eq-g}) for $i = 1, j = 2$ gives
  \begin{eqnarray*}
    h & = & \log (\beta x_1) - \beta x_1 + 1
  \end{eqnarray*}
  which is negative unless $x_1$ also equals $1 / \beta$. Yet, $x_1 = 1 /
  \beta$ implies $\beta = q$ and $h = 0$. But $\tmmathbf{x}= (1 / q, \ldots, 1
  / q)$ is a minimizer of the free energy $f_{\beta, 0}$ only for $\beta
  \leqslant \beta_c$. As $q > 2 \Rightarrow \beta_c < q$, the only case of
  equality is $q = \beta = 2$ and $h = 0$.
\end{proof}

In a second Lemma we compare the counting measure on $\Omega^n$ with the
Lebesgue measure. This will help in the proofs of Theorems \ref{thm-fluct1},
\ref{thm-fluct-cond}, \ref{thm-probacoex} and \ref{thm-fluct-b0h0}. We denote
by $\mathcal{L}$ the Lebesgue measure on hyperplanes.

\begin{lemma}
  \label{lem-Leb}Let $\Pi : \mathbbm{R}^q \rightarrow \mathbbm{R}^q$ be a
  affine and one-to-one transformation. Let $P = [0, 1]^d \cap \mathcal{H}$.
  Then, for any $f : \mathbbm{R}^q \rightarrow \mathbbm{R}$ bounded,
  \begin{eqnarray*}
    \sum_{\tmmathbf{X} \in \Omega^n} f (\tmmathbf{X}) & \leqslant &
    \frac{1}{\mathcal{L}(\Pi (n^{- 1} P))}  \int_{\Pi (\Omega + n^{- 1} P)}
    \varphi \mathd \mathcal{L}
  \end{eqnarray*}
  where
  \begin{eqnarray*}
    \varphi (\tmmathbf{z}) & = & \sup_{\Pi^{- 1} (\tmmathbf{z}) - n^{- 1} P} f
    \text{, \ \ \ \ } \forall \tmmathbf{z} \in \mathbbm{R}^q .
  \end{eqnarray*}
\end{lemma}

\begin{remark}
  Applying this to $- f$ one obtains a useful lower bound.
\end{remark}

\begin{proof}
  For any $\tmmathbf{z} \in \Pi \left( \tmmathbf{x}+ n^{- 1} P \right)$, one
  has $\Pi^{- 1} (\tmmathbf{z}) \in \tmmathbf{x}+ n^{- 1} P$ hence
  $\tmmathbf{x} \in \Pi^{- 1} (\tmmathbf{z}) - n^{- 1} P$. Thus $\varphi
  (\tmmathbf{z}) \geqslant f (\tmmathbf{x})$, and
  \begin{eqnarray*}
    f (\tmmathbf{x}) & \leqslant & \frac{1}{\mathcal{L}(\Pi (n^{- 1} P))}
    \int_{\Pi \left( \tmmathbf{x}+ n^{- 1} P \right)} \varphi d\mathcal{L}.
  \end{eqnarray*}
  The claim follows when we sum over $\tmmathbf{x} \in \Omega^n$, as $\Omega^n
  + n^{- 1} P = \Omega + n^{- 1} P$.
\end{proof}

\subsection{Gaussian fluctuations}

The limit theorems will be proved as consequences of a Taylor expansion of the
free energy. First we consider a second order expansion of $f_{\beta, h}$,
that will be enough to describe the Gaussian fluctuations at $(\beta, h) \neq
(\beta_0, h_0)$.

This section is organized as follows. First we give a series of Lemmas that
permit to establish Proposition \ref{prop-f2} below. Then we give the proofs
of Theorems \ref{thm-fluct1}, \ref{thm-fluct-cond} and \ref{thm-probacoex}.

\subsubsection{Taylor expansion of the free energy}

The Taylor-Lagrange formula applied to the $\mathcal{C}^{\infty}$ function $t
\in [0, 1] \mapsto f_{\beta, h} (\tmmathbf{x}+ t\tmmathbf{w})$ yields:

\begin{lemma}
  \label{lem-Taylor-x}Let $\tmmathbf{x} \in \Omega$ be a global minimizer of
  $f_{\beta, h}$ and $\tmmathbf{w} \in \mathcal{H}$ such that
  $\tmmathbf{x}+\tmmathbf{w} \in \Omega_{0^+}$. Then, there exists $\alpha \in
  (0, 1)$ such that
  \begin{eqnarray}
    f_{\beta, h} (\tmmathbf{x}+\tmmathbf{w}) & = & f_{\beta, h} (\tmmathbf{x})
    + \frac{1}{2} \sum_{i = 1}^q \left( \frac{1}{x_i + \alpha w_i} - \beta
    \right) w_i^2 . \label{eq-T2} 
  \end{eqnarray}
\end{lemma}

On the other hand, the influence of $\beta_n$ and $h_n$ is immediate to
characterize:

\begin{lemma}
  \label{lem-Taylor-bh}For any $\beta, \beta_n, h, h_n$ and any $\tmmathbf{x}
  \in \Omega$, the following equality holds:
  \begin{eqnarray}
    f_{\beta_n, h_n} (\tmmathbf{x}) & = & f_{\beta, h} (\tmmathbf{x}) -
    \frac{\beta_n - \beta}{2} \sum_{i = 1}^q x_i^2 - (h_n - h) x_1 . 
  \end{eqnarray}
\end{lemma}

\subsubsection{The quadratic form}

Given $\tmmathbf{x} \in \Omega_{0^+}$ and $\beta \geqslant 0$ we consider the
quadratic form $Q_{\tmmathbf{x}, \beta} : \mathcal{H} \mapsto \mathbbm{R}$
defined by
\begin{eqnarray}
  Q_{\tmmathbf{x}, \beta} (\tmmathbf{w}) & = & \sum_{i = 1}^q \left(
  \frac{1}{x_i} - \beta \right) w_i^2 .  \label{defQ}
\end{eqnarray}
This is the quadratic form that appears in Lemma \ref{lem-Taylor-x}. When it
is positive definite it determines the fluctuations. We have:

\begin{lemma}
  \label{lem-quad-form}Let $\beta, h \geqslant 0$ and $\tmmathbf{x}$ be some
  global minimizer of $f_{\beta, h}$.
  \begin{enumerateroman}
    \item The quadratic form $Q_{\tmmathbf{x}, \beta}$ is positive definite on
    $\mathcal{H}$ if and only $(\beta, h) \neq (\beta_0, h_0)$.
    
    \item When $(\beta, h) = (\beta_0, h_0)$ and
    $\tmmathbf{x}=\tmmathbf{x}_0$, the kernel of $Q_{\tmmathbf{x}, \beta}$ is
    $\tmop{Vect} (\tmmathbf{u})$ where
    \[ \tmmathbf{u}= \left( 1 - q, 1, \ldots, 1 \right) . \]
  \end{enumerateroman}
\end{lemma}

\begin{proof}
  First we assume that $(\beta, h) \neq (\beta_0, h_0)$ and prove that
  $Q_{\tmmathbf{x}, \beta}$ is positive definite. According to Proposition
  \ref{prop-struct-x} the vector $\tmmathbf{x}$ has one coordinate repeated at
  least $q - 1$ times. Let $j \in \{1, \ldots, q\}$ be the smallest index such
  that $x_j = \max (x_i)$, and $J =\{1, \ldots, q\} \setminus \{j\}$. For any
  $\tmmathbf{w} \in \mathcal{H}$ one has $w_j = - \sum_{i \in J} w_i$, hence
  \begin{eqnarray*}
    Q_{\tmmathbf{x}, \beta} (\tmmathbf{w}) & = & \left( \frac{1}{\max (x_i)} -
    \beta \right) \left( \sum_{i \in J} w_i \right)^2 + \left( \frac{1}{\min
    (x_i)} - \beta \right) \sum_{i \in J} w_i^2 .
  \end{eqnarray*}
  Now we let
  \begin{eqnarray*}
    \alpha_j (\tmmathbf{w}) & = & \frac{1}{q - 1}  \frac{\left( \sum_{i \in I}
    w_i \right)^2}{\sum_{i \in J} w_i^2},
  \end{eqnarray*}
  which belongs to the interval $[0, 1]$ according to Cauchy-Schwarz
  inequality, and obtain
  \begin{eqnarray}
    Q_{\tmmathbf{x}} (\tmmathbf{w}) & = & \left[ \left( \frac{1}{\max (x_i)} -
    \beta \right) (q - 1) \alpha_j (\tmmathbf{w}) + \left( \frac{1}{\min
    (x_i)} - \beta \right) \right] \sum_{i \in J} w_i^2 .  \label{quadform}
  \end{eqnarray}
  Hence the quadratic form $Q_{\tmmathbf{x}}$ is positive definite on
  $\mathcal{H}$ if and only if the factor in (\ref{quadform}) is strictly
  positive at both $\alpha = 0$ and $\alpha = 1$, that is to say if
  \begin{eqnarray}
    \frac{1}{\min (x_i)} - \beta & > & 0  \label{condD2fposa}\\
    \text{and \ } \frac{q - 1}{\max (x_i)} + \frac{1}{\min (x_i)} - q \beta &
    > & 0.  \label{condD2fposb}
  \end{eqnarray}
  Condition (\ref{condD2fposa}) is true as $(\beta, h) \neq (\beta_0, h_0)$,
  cf. Proposition \ref{prop-struct-x}. Condition (\ref{condD2fposb}) is
  equivalent to
  \begin{eqnarray*}
    1 - q \beta \min (x_i) \max (x_i) & > & 0
  \end{eqnarray*}
  as $(q - 1) \min (x_i) + \max (x_i) = 1$. The reader will remark that $\min
  (x_i) \max (x_i)$ is constant over all the minimizers of $f_{\beta, h}$
  described in Theorem \ref{thm-transition}. Hence we might take the minimizer
  of the form $\tmmathbf{x}=\tmmathbf{x}_z$ as in (\ref{xz}), that is $x_1 =
  (1 + z) / 2$ and $x_2 = \cdots = x_q = (1 - z) / (q - 1)$, which reveals
  that condition (\ref{condD2fposb}) is equivalent to $z \mapsto f_{\beta, h}
  (\tmmathbf{x}_z)$ having a positive second derivative at its minima, as
  \[ \frac{\mathd^2 f_{\beta, h} (\tmmathbf{x}_z)}{\mathd z^2} = \frac{1}{1 -
     z^2} - \frac{\beta q}{4 (q - 1)}, \]
  which is the case again as $(\beta, h) \neq (\beta_0, h_0)$ (see the
  discussion after (\ref{d2fz})).
  
  Assume now that $(\beta, h) = (\beta_0, h_0)$. If $q = 2$, the quadratic
  form $Q_{\tmmathbf{x}, \beta}$ is identically zero on $\mathcal{H}=
  \tmop{Vect} (\tmmathbf{u})$. If $q \geqslant 3$, we have $h = h_0 > 0$ hence
  $j = 1$. The quadratic form vanishes at $\tmmathbf{w} \in \mathcal{H}$ if
  and only if $\alpha_1 (\tmmathbf{w}) = 1$. In view of the definition of
  $\alpha$, this is the case of equality in the Cauchy-Schwarz inequality:
  $\alpha_1 (\tmmathbf{w}) = 1 \Leftrightarrow w_2 = \cdots = w_q
  \Leftrightarrow \tmmathbf{w} \in \tmop{Vect} (\tmmathbf{u})$.
\end{proof}

\subsubsection{Centering of the fluctuations}

As in Theorems \ref{thm-fluct1} and \ref{thm-fluct-cond} we let
$\tmmathbf{d}_n$ the smallest $\tmmathbf{d} \in \mathcal{H}$ such that
$\tmmathbf{x}+\tmmathbf{d} \in \Omega$ is a global minimizer of $f_{\beta_n,
h_n}$. We have:

\begin{proposition}
  \label{prop-f2}Assume that $(\beta_n, h_n) \rightarrow (\beta, h)$ and let
  $\tmmathbf{x} \in \Omega$ be a global minimizer of $f_{\beta, h}$.
  \begin{enumerateroman}
    \item For any $R > 0$,
    \begin{eqnarray}
      nf_{\beta_n, h_n} \left( \tmmathbf{x}+\tmmathbf{d}_n + n^{- 1 / 2}
      \tmmathbf{w} \right) & = & nf_{\beta_n, h_n} \left(
      \tmmathbf{x}+\tmmathbf{d}_n \right) + \frac{1}{2}
      \tmmathbf{Q}_{\tmmathbf{x}, \beta} (\tmmathbf{w}) + o_n (1) 
      \label{eq-f2dev}
    \end{eqnarray}
    uniformly over $\tmmathbf{w} \in \mathcal{H} \cap B (0, R)$.
    
    \item If $(\beta, h) \neq (\beta_0, h_0)$, for small enough $\varepsilon >
    0$ there is $\lambda > 0$ such that, for $n$ large enough and any
    $\tmmathbf{w} \in \mathcal{H}$ with $\|\tmmathbf{w}\| \leqslant
    \varepsilon n^{1 / 2}$,
    \begin{eqnarray}
      nf_{\beta_n, h_n} \left( \tmmathbf{x}+\tmmathbf{d}_n + n^{- 1 / 2}
      \tmmathbf{w} \right) & \geqslant & nf_{\beta_n, h_n} \left(
      \tmmathbf{x}+\tmmathbf{d}_n \right) + \lambda \|\tmmathbf{w}\|^2 . 
      \label{eq-f2tight}
    \end{eqnarray}
  \end{enumerateroman}
\end{proposition}

\begin{proof}
  We begin with an application of Lemma \ref{lem-Taylor-x} at the global
  minimum point $\tmmathbf{x}+\tmmathbf{d}_n$:
  \begin{eqnarray*}
    nf_{\beta_n, h_n} \left( \tmmathbf{x}+\tmmathbf{d}_n + n^{- 1 / 2}
    \tmmathbf{w} \right) & = & nf_{\beta_n, h_n} \left(
    \tmmathbf{x}+\tmmathbf{d}_n \right) + \frac{1}{2}
    \tmmathbf{Q}_{\tmmathbf{x}+\tmmathbf{d}_n + \alpha n^{- 1 / 2}
    \tmmathbf{w}, \beta} (\tmmathbf{w})
  \end{eqnarray*}
  for some $\alpha \in (0, 1)$ depending on $n$ and $\tmmathbf{w}$. For
  (\ref{eq-f2dev}) we only have to notice that $\tmmathbf{d}_n + \alpha n^{- 1
  / 2} \tmmathbf{w}= o (1)$. For (\ref{eq-f2tight}) we remark that as
  $\tmmathbf{d}_n \rightarrow 0$ (cf. Remark \ref{rmk-d}), for all $n$ large
  enough and $\|\tmmathbf{w}\| \leqslant \varepsilon n^{1 / 2}$,
  \begin{eqnarray*}
    \left\| \tmmathbf{d}_n + \alpha n^{- 1 / 2} \tmmathbf{w} \right\| &
    \leqslant & 2 \varepsilon .
  \end{eqnarray*}
  As this can be made arbitrary small, for small enough $\varepsilon$ the
  quadratic form $\tmmathbf{Q}_{\tmmathbf{x}+\tmmathbf{d}_n + \alpha n^{- 1 /
  2} \tmmathbf{w}, \beta}$ dominates $\tmmathbf{Q}_{x, \beta} / 2$, which is
  definite positive after Lemma \ref{lem-quad-form}.
\end{proof}

\subsubsection{Some linear algebra}

The next Lemma will be useful at the time of computing inverses or
determinants. Denote by $I_n$ the $n \times n$ unitary matrix and by $A_n$ the
$n \times n$ matrix with all entries equal to $1$.

\begin{lemma}
  \label{lem-ApI}Let $M = aA_n + bI_n$.
  \begin{enumerateroman}
    \item The determinant of $M$ is
    \[ \det \left( M \right) = b^{n - 1}  \left( b + na \right) . \]
    \item When $M$ is invertible, it has
    \[ M^{- 1} = \left\{ \begin{array}{ll}
         (a + b)^{- 1} & \text{if $n = 1$}\\
         \frac{1}{b} \left( I_n - \frac{a}{na + b} A_n \right) & \text{if } n
         \geqslant 2.
       \end{array} \right. \]
  \end{enumerateroman}
\end{lemma}

\begin{proof}
  We prove the first point as follows: let $P (\lambda) = \det (\lambda I_n -
  A_n)$ be the characteristic polynomial for the matrix $A_n$. The matrix
  $A_n$ has rank $1$ and eigenvalues $0, \ldots, 0, n$, which are the roots of
  the unitary polynomial $P$, thus $P (\lambda) = \lambda^{n - 1} (\lambda -
  n)$. The second point follows from an immediate computation.
\end{proof}

\subsubsection{Proof of Theorems \ref{thm-fluct1}, \ref{thm-fluct-cond} and
\ref{thm-probacoex}}

As a consequence of Proposition \ref{prop-f2} we give the proof of Theorems
\ref{thm-fluct1}, \ref{thm-fluct-cond} and \ref{thm-probacoex}.

\begin{proof}
  (Theorem \ref{thm-fluct1}). First we condition $\tmmathbf{W}= n^{- 1 / 2}
  (\tmmathbf{N}- n\tmmathbf{x}- n\tmmathbf{d}_n)$ on $\|\tmmathbf{W}\|< R$ for
  some positive $R$. For $g : \mathbbm{R}^q \mapsto \mathbbm{R}$ continuous
  bounded, Lemma \ref{lem-stirling} and (\ref{eq-f2dev}) in Proposition
  \ref{prop-f2} yield
  \begin{eqnarray*}
    Z_{\beta_n, h_n, n} \mu_{\beta_n, h_n, n} \left( g
    (\tmmathbf{W})\tmmathbf{1}_{\{\|\tmmathbf{W}\| \leqslant R\}} \right)  & =
    & (1 + o_n (1)) n^{- \frac{q - 1}{2}} A_{\beta} (\tmmathbf{x}) e^{-
    nf_{\beta_n, h_n} \left( \tmmathbf{x}+\tmmathbf{d}_n \right)}\\
    &  & \times \sum_{\tmmathbf{N}/ n \in \Omega^n} g
    (\tmmathbf{W})\tmmathbf{1}_{\{\|\tmmathbf{W}\| \leqslant R\}} e^{-
    \frac{1}{2} \tmmathbf{Q}_{\tmmathbf{x}, \beta} (\tmmathbf{W})} .
  \end{eqnarray*}
  The transformation $\tmmathbf{X}=\tmmathbf{N}/ n \mapsto \Pi (\tmmathbf{X})
  =\tmmathbf{W}$ is affine. The image of $\Omega$ by $\Pi$ is greater than
  $\mathcal{H} \cap B (0, R)$ for large enough $n$ as $\tmmathbf{x} \in
  \Omega_{0^+}$, and on the other hand $\Pi (\tmmathbf{X}+ n^{- 1} P) = \Pi
  (\tmmathbf{X}) + n^{- 1 / 2} P$, that is to say the dimensions of the image
  of the lattice element $P$ go to zero. Hence Lemma \ref{lem-Leb} gives
  \begin{eqnarray}
    Z_{\beta_n, h_n, n} \mu_{\beta_n, h_n, n} \left( g
    (\tmmathbf{W})\tmmathbf{1}_{\{\|\tmmathbf{W}\| \leqslant R\}} \right)  & =
    & (1 + o_n (1)) A_{\beta} (\tmmathbf{x}) e^{- nf_{\beta_n, h_n} \left(
    \tmmathbf{x}+\tmmathbf{d}_n \right)} \times \nonumber\\
    &  & \frac{1}{\mathcal{L}(P)} \int_{\mathcal{H} \cap B (0, R)} g
    (\tmmathbf{w}) e^{- \frac{1}{2} \tmmathbf{Q}_{\tmmathbf{x}, \beta}
    (\tmmathbf{w})} \mathd \mathcal{L}(\tmmathbf{w}) .  \label{equiv-Zmug}
  \end{eqnarray}
  In other words the law of $\tmmathbf{W}$ conditioned on $\|\tmmathbf{W}\|
  \leqslant R$ converges to the distribution on $\mathcal{H} \cap B (0, R)$
  with density proportional to
  \begin{equation}
    \tmmathbf{w} \mapsto e^{- \frac{1}{2} Q_{\tmmathbf{x}, \beta}
    (\tmmathbf{w})} \label{eq-densQ}
  \end{equation}
  with respect to the Lebesgue measure on $\mathcal{H} \cap B (0, R)$.

  Now we show that the variable $\tmmathbf{W}$ is tight. Thanks to Theorem
  \ref{thm-LDP} we know already that for any $\varepsilon > 0$,
  \begin{eqnarray*}
    \limsup_n \mu_{\beta_n, h_n, n} \left( \left\| \tmmathbf{N}- n\tmmathbf{x}
    \right\| \geqslant \varepsilon n \right) & = & 0.
  \end{eqnarray*}
  Thus it is enough to show that, for small enough $\varepsilon > 0$,
  \begin{eqnarray}
    \lim_{\kappa \rightarrow \infty} \limsup_n \mu_{\beta_n, h_n, n} \left(
    \|\tmmathbf{W}\| \geqslant \kappa | \|\tmmathbf{W}\| \leqslant \varepsilon
    n^{1 / 2} \right) & = & 0.  \label{eq-Wtight}
  \end{eqnarray}
  According to Lemma \ref{lem-stirling} and Proposition \ref{prop-f2}, for
  small enough $\varepsilon > 0$ and large enough $\kappa$ there is $\lambda >
  0$ such that
  \begin{eqnarray*}
    \mu_{\beta_n, h_n, n} \left( \|\tmmathbf{W}\| \geqslant \kappa |
    \|\tmmathbf{W}\| \leqslant \varepsilon n^{1 / 2} \right) & \leqslant &
    \frac{\mu_{\beta_n, h_n, n} \left( \kappa \leqslant \|\tmmathbf{W}\|
    \leqslant \varepsilon n^{1 / 2} \right)}{\mu_{\beta_n, h_n, n} \left(
    \|\tmmathbf{W}\| \leqslant \kappa \right)}\\
    & \leqslant & (1 + o_n (1)) \frac{\sum_{\tmmathbf{N}/ n \in \Omega^n :
    \kappa \leqslant \|\tmmathbf{W}\| \leqslant \varepsilon n^{1 / 2}} e^{-
    \lambda \|\tmmathbf{W}\|^2}}{\sum_{\tmmathbf{N}/ n \in \Omega^n :
    \|\tmmathbf{W}\| \leqslant \kappa} e^{- \frac{1}{2}
    \tmmathbf{Q}_{\tmmathbf{x}, \beta} (\tmmathbf{W})}}\\
    & \leqslant & (1 + o_n (1)) \frac{\int_{\mathcal{H} \setminus B (0,
    \kappa)} e^{- \lambda \|\tmmathbf{w}\|^2} \mathd
    \mathcal{L}(\tmmathbf{w})}{\int_{\mathcal{H} \cap B (0, \kappa)} e^{-
    \frac{1}{2} \tmmathbf{Q}_{x, \beta} (\tmmathbf{w})} \mathd
    \mathcal{L}(\tmmathbf{w})} .
  \end{eqnarray*}
  after Lemma \ref{lem-Leb}. Since $\tmmathbf{Q}_{\tmmathbf{x}, \beta}$ is
  positive definite the ratio goes to $0$ as $\kappa \rightarrow \infty$,
  giving (\ref{eq-Wtight}).

  Let us show that this limit distribution is as well the distribution of the
  centered Gaussian vector with covariance matrix (\ref{defKx}). We take
  $\tmmathbf{V}$ a random vector in $\mathcal{H}$ with the density
  (\ref{eq-densQ}) with respect to the Lebesgue measure on $\mathcal{H}$. The
  law of $\tmmathbf{V}$ is also proportional to
  \[ e^{- \frac{1}{2} Q_{\tmmathbf{x}, \beta} (\tmmathbf{v})} \mathd v_2
     \cdots \mathd v_q . \]
  This density can be expressed only in terms of the truncated vector
  $\tilde{\tmmathbf{V}} = (V_2, \ldots, V_q)$. Indeed, if we take
  \begin{eqnarray}
    H & = & \left( \frac{1}{x_1} - \beta \right) A_{q - 1} + \left(
    \frac{1}{x_q} - \beta \right) I_{q - 1}  \label{defH}
  \end{eqnarray}
  we have $Q_{\tmmathbf{x}, \beta} (\tmmathbf{v}) = \text{}^t
  \tilde{\tmmathbf{v}} H \widetilde{\tmmathbf{v}}$ and thus the covariance
  matrix of $\tilde{\tmmathbf{V}}$ is
  \begin{eqnarray*}
    H^{- 1} & = & \left( \frac{1}{x_1} - \beta \right)^{- 1} \times \left(
    I_{q - 1} - \frac{\frac{1}{x_1} - \beta}{\frac{1}{x_1 x_q} - q \beta} A_{q
    - 1} \right)
  \end{eqnarray*}
  according to Lemma \ref{lem-ApI}. Using the relation $V_1 = - \sum_{i = 2}^q
  V_i$ we compute the remaining covariance coefficients, leading to the
  completed matrix (\ref{defKx}). The rank of the matrix is not less that that
  of $H$, that is $q - 1$, and it is also strictly less than $q$ because of
  the linear constraint $\tmmathbf{V} \in \mathcal{H}$.
\end{proof}

\begin{proof}
  (Theorem \ref{thm-fluct-cond}). The former proof can be repeated almost
  verbatim. One has to take care however that $\tmmathbf{x}$ needs not be of
  the particular form $x_2 = \cdots = x_q$ (although it still has a coordinate
  repeated $q - 1$ times), and that the variable which is tight is
  $\tmmathbf{W}$ {\tmem{conditioned}} on $\tmmathbf{N}/ n \in B (\tmmathbf{x},
  \varepsilon)$.
\end{proof}

\begin{proof}
  (Theorem \ref{thm-probacoex}). The tightness of $\tmmathbf{W}$ conditioned
  on $\tmmathbf{N}/ n \in B (\tmmathbf{x}, \varepsilon)$ and the convergence
  of the law of $\tmmathbf{W}$ on bounded sets (cf. the proof of Theorem
  \ref{thm-fluct1}) imply that for any $\varepsilon > 0$ smaller than the
  distance between any two minimizers of $f_{\beta, h}$,
  \begin{eqnarray*}
    \frac{\mu_{\beta_n, h_n, n} \left( \frac{\tmmathbf{N}}{n} \in B
    (\tmmathbf{x}, \varepsilon) \right)}{\mu_{\beta_n, h_n, n} \left(
    \frac{\tmmathbf{N}}{n} \in B (\tmmathbf{x}', \varepsilon) \right)} & = &
    (1 + o_n (1)) \frac{A_{\beta} (\tmmathbf{x}) e^{- nf_{\beta_n, h_n} \left(
    \tmmathbf{x}+\tmmathbf{d}_n (\tmmathbf{x}) \right)}  \int_{\mathcal{H}}
    e^{- \frac{1}{2} \tmmathbf{Q}_{\tmmathbf{x}, \beta} (\tmmathbf{w})} \mathd
    \mathcal{L}(\tmmathbf{w})}{A_{\beta} (\tmmathbf{x}') e^{- nf_{\beta_n,
    h_n} \left( \tmmathbf{x}' +\tmmathbf{d}_n (\tmmathbf{x}') \right)}
    \int_{\mathcal{H}} e^{- \frac{1}{2} \tmmathbf{Q}_{\tmmathbf{x}', \beta}
    (\tmmathbf{w})} \mathd \mathcal{L}(\tmmathbf{w})} .
  \end{eqnarray*}
  Hence we call
  \begin{eqnarray*}
    C_{\beta_n, h_n, n} (\tmmathbf{x}) & = & A_{\beta} (\tmmathbf{x}) e^{-
    nf_{\beta_n, h_n} \left( \tmmathbf{x}+\tmmathbf{d}_n (\tmmathbf{x})
    \right)}  \int_{\mathcal{H}} e^{- \frac{1}{2} \tmmathbf{Q}_{\tmmathbf{x},
    \beta} (\tmmathbf{w})} \mathd \mathcal{L}(\tmmathbf{w})
  \end{eqnarray*}
  and give an equivalent to $C_{\beta_n, h_n, n} (\tmmathbf{x})$. First we
  compute the integral up to a constant factor. We pick $j \in \{1, \ldots,
  q\}$ such that $x_j = \max (x_i)$ and let $J =\{1, \ldots, q\} \setminus
  \{j\}$. The Lebesgue measure on $\mathcal{H}$ is proportional to the measure
  induced on $\tmmathbf{w}$ by $\prod_{i \in J} \mathd w_i$, given $w_j = -
  \sum_{i \in J} w_i$. As in (\ref{defH}) we let
  \begin{eqnarray*}
    H_{\tmmathbf{x}, \beta} & = & \left( \frac{1}{\max (x_i)} - \beta \right)
    A_{q - 1} + \left( \frac{1}{\min (x_i)} - \beta \right) I_{q - 1}
  \end{eqnarray*}
  and $\tilde{\tmmathbf{w}} = (w_i)_{i \in J}$, thus $\text{}^t
  \tilde{\tmmathbf{w}} H_{\tmmathbf{x}, \beta} \widetilde{\tmmathbf{w}} =
  Q_{\tmmathbf{x}, \beta} (\tmmathbf{w})$ and therefore
  \begin{eqnarray*}
    \int_{\mathcal{H}} e^{- \frac{1}{2} \tmmathbf{Q}_{\tmmathbf{x}, \beta}
    (\tmmathbf{w})} \prod_{i \in J} \mathd w_i & = & \sqrt{2 \pi}^{q - 1}
    \sqrt{\det H_{\tmmathbf{x}, \beta}^{- 1}}\\
    & = & \sqrt{2 \pi}^{q - 1} \left[ \left( \frac{1}{\min (x_i)} - \beta
    \right)^{q - 2}  \left( \frac{1}{\min (x_i)} + \frac{q - 1}{\max (x_i)} -
    q \beta \right) \right]^{- 1 / 2}
  \end{eqnarray*}
  according to Lemma \ref{lem-ApI}. If we multiply with the prefactor
  $A_{\beta} (\tmmathbf{x})$ we obtain
  \begin{eqnarray}
    A_{\beta} (\tmmathbf{x}) \int_{\mathcal{H}} e^{- \frac{1}{2}
    \tmmathbf{Q}_{\tmmathbf{x}, \beta} (\tmmathbf{w})} \prod_{i \in J} \mathd
    w_i & = & e^{- \frac{\beta}{2}} (1 - \beta \min (x_i))^{\frac{2 - q}{2}} 
    \left( 1 - q \beta \max (x_i) \min (x_i) \right)^{- 1 / 2} 
    \label{eq-Aint}
  \end{eqnarray}
  as $\max (x_i) + (q - 1) \min (x_i) = 1$. Then we use Lemma
  \ref{lem-Taylor-bh}:
  \begin{eqnarray*}
    nf_{\beta_n, h_n} \left( \tmmathbf{x}+\tmmathbf{d}_n (\tmmathbf{x})
    \right) & = & nf_{\beta, h} \left( \tmmathbf{x}+\tmmathbf{d}_n
    (\tmmathbf{x}) \right) - n \frac{\beta_n - \beta}{2} \sum_{i = 1}^q \left(
    x_i + d_{n, i} (\tmmathbf{x}) \right)^2 - n (h_n - h) \left( x_1 + d_1
    (\tmmathbf{x}) \right)\\
    & = & nf_{\beta, h} \left( \tmmathbf{x} \right) - \frac{\lambda}{2}
    \sum_{i = 1}^q x_i^2 - \nu x_1 + o_n (1)
  \end{eqnarray*}
  as $\tmmathbf{d}_n (\tmmathbf{x}) = O (1 / n)$ and $\tmmathbf{x}$ is a
  global minimizer of $f_{\beta, h}$. Thus we have shown that
  \begin{eqnarray*}
    C_{\beta_n, h_n, n} (\tmmathbf{x}) & = & (1 + o_n (1)) \frac{\mathd
    \mathcal{L}}{\prod_{i \in J} \mathd w_i} \exp \left( - nf_{\beta, h}
    \left( \tmmathbf{x} \right) + \frac{\lambda}{2} \sum_{i = 1}^q x_i^2 + \nu
    x_1 - \frac{\beta}{2} \right) \\
    &  & \times (1 - \beta \min (x_i))^{\frac{2 - q}{2}}  \left( 1 - q \beta
    \max (x_i) \min (x_i) \right)^{- 1 / 2}
  \end{eqnarray*}
  where the factor $\mathd \mathcal{L}/ \prod_{i \in J} \mathd w_i$ does not
  depend on $J$. The claim follows from the remark that the product $\min
  (x_i) \max (x_i)$ is constant over all the global minimizers $\tmmathbf{x},
  \tmmathbf{x}', \ldots$ of the free energy $f_{\beta, h}$ at any $(\beta, h)
  \in h_T$, cf. Theorem \ref{thm-transition}.
\end{proof}

\subsection{Limit theorems at criticality}

The proof of Theorem \ref{thm-fluct-b0h0} relies again on a Taylor expansion
of the free energy:

\begin{lemma}
  \label{lem-T4}Let $(\beta, h) = (\beta_0, h_0)$ and
  $\tmmathbf{x}=\tmmathbf{x}_0 = (1 / 2, 1 / 2 (q - 1), \ldots, 1 / 2 (q -
  1))$ be the unique minimizer of $f_{\beta, h}$. Let $\tmmathbf{u}= (1 - q,
  1, \ldots, 1)$. For all $t \in \mathbbm{R}$ and $\tmmathbf{v} \in
  \mathcal{H} \cap \tmmathbf{u}^{\perp}$ such that $\tmmathbf{x}+
  t\tmmathbf{u}+\tmmathbf{v} \in \Omega_{0^+}$, there are $\alpha, \alpha' \in
  (0, 1)$ such that
  \begin{eqnarray*}
    f_{\beta, h} (\tmmathbf{x}+ t\tmmathbf{u}+\tmmathbf{v}) & = & f_{\beta, h}
    (\tmmathbf{x}) + \frac{1}{2} Q_{\tmmathbf{x}+ t\tmmathbf{u}+ \alpha
    \tmmathbf{v}, \beta} \left( \tmmathbf{v} \right) + \frac{t^4}{12} \sum_{i
    = 1}^q \frac{u_i^4}{(x_i + \alpha' tu_i)^3} .
  \end{eqnarray*}
  Furthermore,
  \begin{eqnarray*}
    \frac{1}{2} Q_{\tmmathbf{x}, \beta} \left( \tmmathbf{v} \right) & = &
    \frac{(q - 1) (q - 2)}{q} \|\tmmathbf{v}\|^2\\
    \text{and \ } \frac{1}{12} \sum_{i = 1}^q \frac{u_i^4}{x_i^3} & = &
    \frac{4}{3} (q - 1)^3
  \end{eqnarray*}
\end{lemma}

\begin{proof}
  A second-order Taylor expansion in $\tmmathbf{v}$ yields
  \begin{eqnarray*}
    f_{\beta, h} (\tmmathbf{x}+ t\tmmathbf{u}+\tmmathbf{v}) & = & f_{\beta, h}
    (\tmmathbf{x}+ t\tmmathbf{u}) + \nabla f_{\beta, h} (\tmmathbf{x}+
    t\tmmathbf{u}) \cdot \tmmathbf{v}+ \frac{1}{2} Q_{\tmmathbf{x}+
    t\tmmathbf{u}+ \alpha \tmmathbf{v}, \beta} \left( \tmmathbf{v} \right)
  \end{eqnarray*}
  for some $\alpha \in (0, 1)$. The last $q - 1$ coordinates of the gradient
  $\nabla f_{\beta, h} (\tmmathbf{x}+ t\tmmathbf{u})$ are equal, hence it is
  orthogonal to $\tmmathbf{v}$. Then a fourth order expansion in $t$ gives
  \begin{eqnarray}
    f_{\beta, h} (\tmmathbf{x}+ t\tmmathbf{u}) & = & f_{\beta, h}
    (\tmmathbf{x}) + \frac{t^4}{12} \sum_{i = 1}^q \frac{u_i^4}{(x_i + \alpha'
    tu_i)^3}  \label{eq-f4}
  \end{eqnarray}
  for some $\alpha' \in (0, 1)$. Indeed, the first order term is zero as
  $\tmmathbf{x}$ is the global minimizer of $f_{\beta, h}$. The second order
  term is $Q_{\tmmathbf{x}, \beta} \left( t\tmmathbf{u} \right) / 2 = 0$ in
  view of Lemma \ref{lem-quad-form}. Hence the third order term is $0$,
  yielding (\ref{eq-f4}).
  
  Let us prove the last two formulas. The assumption $\tmmathbf{v} \in
  \mathcal{H} \cap \tmmathbf{u}^{\perp}$ implies $v_1 = 0$, hence
  \begin{eqnarray*}
    \frac{1}{2} Q_{\tmmathbf{x}, \beta} \left( \tmmathbf{v} \right) & = &
    \frac{1}{2} \sum_{i = 2}^q \left( \frac{1}{x_i} - \beta \right) v_i^2\\
    & = & \frac{1}{2} \sum_{i = 2}^q \left( 2 (q - 1) - 4 \frac{q - 1}{q}
    \right) v_i^2\\
    & = & \frac{(q - 1) (q - 2)}{q}  \sum_{i = 1}^q v_i^2 .
  \end{eqnarray*}
  On the other hand:
  \[ \sum_{i = 1}^q \frac{u_i^4}{x_i^3} = 8 (q - 1)^4 + (q - 1) 8 (q - 1)^3 =
     16 (q - 1)^4 . \]
\end{proof}

Using Lemma \ref{lem-T4} we establish the analog of Proposition \ref{prop-f2}:

\begin{proposition}
  \label{prop-f4dev}Assume that $(\beta_n, h_n) \rightarrow (\beta_0, h_0)$
  with $\beta_n - \beta_0 = o (n^{- 3 / 4})$ and $h_n - h_0 = o (n^{- 3 /
  4})$, and let $\tmmathbf{x}_0 \in \Omega$ be the unique global minimizer of
  $f_{\beta_0, h_0}$.
  \begin{enumerateroman}
    \item For any $R > 0$,
    \begin{eqnarray}
      nf_{\beta_n, h_n} \left( \tmmathbf{x}+ n^{- 1 / 4} t\tmmathbf{u}+ n^{- 1
      / 2} \tmmathbf{v} \right) & = & nf_{\beta, h} \left( \tmmathbf{x}
      \right) - n \frac{\beta_n - \beta_0}{2} \|\tmmathbf{x}\|^2 - n (h_n -
      h_0) x_1 + \nonumber\\
      &  & \frac{(q - 1) (q - 2)}{q} \|\tmmathbf{v}\|^2 + \frac{4}{3} (q -
      1)^3 t^4 + o_n (1)  \label{eq-f4dev}
    \end{eqnarray}
    uniformly over $\tmmathbf{v} \in \mathcal{H} \cap \tmmathbf{u}^{\perp}
    \cap B (0, R)$ and $t \in [- R, R]$.
    
    \item For small enough $\varepsilon > 0$ and large enough $R$, for $n$
    large enough, for any $\tmmathbf{v} \in \mathcal{H} \cap
    \tmmathbf{u}^{\perp}$, $t \in \mathbbm{R} \setminus [- R, R]$ such that
    $\|n^{- 1 / 4} t\tmmathbf{u}+ n^{- 1 / 2} \tmmathbf{v}\| \leqslant
    \varepsilon$,
    \begin{eqnarray}
      nf_{\beta_n, h_n} \left( \tmmathbf{x}+ n^{- 1 / 4} t\tmmathbf{u}+ n^{- 1
      / 2} \tmmathbf{v} \right) & \geqslant & nf_{\beta, h} \left(
      \tmmathbf{x} \right) - n \frac{\beta_n - \beta_0}{2} \|\tmmathbf{x}\|^2
      - n (h_n - h_0) x_1 \nonumber\\
      &  & +\|\tmmathbf{v}\|^2 / 4 + t^4 / 2  \label{eq-f4tight}
    \end{eqnarray}
  \end{enumerateroman}
\end{proposition}

\begin{proof}
  We first apply Lemma \ref{lem-Taylor-bh}: for $\tmmathbf{w}= n^{- 1 / 4}
  t\tmmathbf{u}+ n^{- 1 / 2} \tmmathbf{v}$, we have
  \begin{eqnarray*}
    nf_{\beta_n, h_n} \left( \tmmathbf{x}+\tmmathbf{w} \right) & = &
    nf_{\beta, h} \left( \tmmathbf{x}+\tmmathbf{w} \right) - n \frac{\beta_n -
    \beta_0}{2} \|\tmmathbf{x}\|^2 - n (h_n - h_0) x_1 + o_n (1)
  \end{eqnarray*}
  uniformly over $\tmmathbf{v} \in \mathcal{H} \cap \tmmathbf{u}^{\perp} \cap
  B (0, R)$ and $t \in [- R, R]$ as $\beta_n - \beta_0 = o (n^{- 3 / 4})$ and
  $h_n - h_0 = o (n^{- 3 / 4})$. Then, Lemma \ref{lem-T4} yields
  \begin{eqnarray*}
    nf_{\beta, h} \left( \tmmathbf{x}+\tmmathbf{w} \right) & = & nf_{\beta, h}
    (\tmmathbf{x}) + \frac{(q - 1) (q - 2)}{q} \|\tmmathbf{v}\|^2 +
    \frac{4}{3} (q - 1)^3 t^4 + o_n (1)
  \end{eqnarray*}
  uniformly over the same domain, and (\ref{eq-f4dev}) follows.
  
  Now we only assume that $\|\tmmathbf{w}\| \leqslant \varepsilon$. For small
  enough $\varepsilon > 0$ we have, after Lemma~\ref{lem-T4}, the lower bound
  \begin{eqnarray}
    nf_{\beta, h} \left( \tmmathbf{x}+\tmmathbf{w} \right) & \geqslant &
    nf_{\beta, h} (\tmmathbf{x}) + \frac{1}{2} \|\tmmathbf{v}\|^2 + t^4 
    \label{eq-f4t1}
  \end{eqnarray}
  (note that, for $q = 2$, $\tmmathbf{v}$ is necessarily $0$). Combining with
  Lemma \ref{lem-Taylor-bh} we obtain that, whenever (\ref{eq-f4t1}) holds,
  \begin{eqnarray*}
    nf_{\beta_n, h_n} \left( \tmmathbf{x}+\tmmathbf{w} \right) & \geqslant &
    c_n (\tmmathbf{x}) +\|\tmmathbf{v}\|^2 \left( \frac{1}{2} - \frac{\beta_n
    - \beta_0}{2} \right) - l_n t - m_n t^2 + t^4
  \end{eqnarray*}
  where
  \begin{eqnarray*}
    c_n (\tmmathbf{x}) & = & nf_{\beta, h} (\tmmathbf{x}) - n \frac{\beta_n -
    \beta_0}{2} \|\tmmathbf{x}\|^2 - n (h_n - h_0) x_1,\\
    l_n & = & n^{3 / 4} (\beta_n - \beta_0)\tmmathbf{x} \cdot \tmmathbf{u}-
    n^{3 / 4} (h_n - h_0) (q - 1)\\
    m_n & = & n^{1 / 2}  \frac{\beta_n - \beta_0}{2} \|\tmmathbf{u}\|^2
  \end{eqnarray*}
  as $\tmmathbf{x}, \tmmathbf{u} \perp \tmmathbf{v}$ and $v_1 = 0$. Now we
  conclude: for any large $n$,
  \begin{eqnarray*}
    \|\tmmathbf{v}\|^2 \left( \frac{1}{2} - \frac{\beta_n - \beta_0}{2}
    \right) & \geqslant & \frac{\|\tmmathbf{v}\|^2}{4} .
  \end{eqnarray*}
  Similarly, as $l_n = o_n (1)$ and $m_n = o (n^{- 1 / 4}) = o_n (1)$, for any
  large $n$ and $t$ large,
  \begin{eqnarray*}
    - l_n t - m_n t^2 + t^4 & \geqslant & \frac{t^4}{2} .
  \end{eqnarray*}
\end{proof}

Finally we give the proof of Theorem \ref{thm-fluct-b0h0}:

\begin{proof}
  (Theorem \ref{thm-fluct-b0h0}). Here we define $\Pi$ as the affine
  transformation such that
  \[ \Pi \left( \tmmathbf{x}+ n^{- 1 / 4} T\tmmathbf{u}+ n^{- 1 / 2}
     \tmmathbf{V}+\tmmathbf{z} \right) =
     T\tmmathbf{u}+\tmmathbf{V}+\tmmathbf{z} \]
  for any $T \in \mathbbm{R}, \tmmathbf{V} \in \mathcal{H} \cap
  \tmmathbf{u}^{\perp}$ and $\tmmathbf{z} \in \mathcal{H}^{\perp}$. It is a
  consequence of Lemmas \ref{lem-stirling}, \ref{lem-Leb} and (\ref{eq-f4dev})
  in Proposition \ref{prop-f4dev} that, conditionally on $\tmmathbf{Z}=
  T\tmmathbf{u}+\tmmathbf{V} \in B (0, R)$, the variable $\tmmathbf{Z}$
  converges in law towards the probability measure on $\mathcal{H} \cap B (0,
  R)$ with density proportional to
  \begin{equation}
    e^{- \frac{(q - 1) (q - 2)}{q} \|\tmmathbf{v}\|^2 - \frac{4 (q - 1)^4}{3}
    t^4} \label{eq-dens4}
  \end{equation}
  with respect to the Lebesgue measure on $\mathcal{H} \cap B (0, R)$. Here
  again, the variable $\tmmathbf{Z}$ is tight thanks to Lemmas
  \ref{lem-stirling}, \ref{lem-Leb} and Proposition \ref{prop-f4dev}.

  The probability measure on $t\tmmathbf{u}+\tmmathbf{v} \in \mathcal{H}$,
  $\tmmathbf{v} \perp \tmmathbf{u}$ with density (\ref{eq-dens4}) has a simple
  structure. It is clear that $T$ and $\tmmathbf{V}$ are independent. The
  vector $\tmmathbf{V}$ is determined by
  \[ \tilde{\tmmathbf{V}} = (V_3, \ldots, V_q) \]
  which has a density proportional to
  \begin{eqnarray*}
    e^{- \frac{(q - 1) (q - 2)}{q} \| \tilde{\tmmathbf{v}} \|^2} \mathd v_3
    \cdots \mathd v_q & = & e^{- \frac{1}{2}  \text{}^t \tilde{\tmmathbf{v}} H
    \tilde{\tmmathbf{v}}} \mathd v_3 \cdots \mathd v_q
  \end{eqnarray*}
  where
  \begin{eqnarray*}
    H & = & 2 \frac{(q - 1) (q - 2)}{q}  \left( A_{q - 2} + I_{q - 2} \right)
    .
  \end{eqnarray*}
  Thus $\tilde{\tmmathbf{V}}$ is the centered Gaussian vector with covariance
  matrix
  \begin{eqnarray*}
    H^{- 1} & = & \frac{q}{2 (q - 1) (q - 2)}  \left( I_{q - 2} - \frac{1}{q -
    1} A_{q - 2} \right) .
  \end{eqnarray*}
  The covariance matrix for $\tmmathbf{V}$ is computed according to $V_1 = 0$
  and $V_2 = - \sum_{i = 3}^q V_i$.
\end{proof}

\subsection{Consequences on the random-cluster model.}

Here we give the proofs of Corollary \ref{cor-proba-gc} and Proposition
\ref{prop-Zrc}:

\begin{proof}
  (Corollary \ref{cor-proba-gc}). When $q > 2$ is an integer, at the critical
  point $(\beta, h) = (\beta_c, 0)$ there are $q + 1$ minimizers for the free
  energy $f_{\beta, h}$, which are, on the one hand, the symmetric state
  \[ \tmmathbf{x}^s = \left( \frac{1}{q}, \cdots, \frac{1}{q} \right) \]
  and on the other hand, the $q$ permutation $\tmmathbf{x}^{a, i}$ of the
  asymmetric state
  \[ \tmmathbf{x}^{a, 1} = \left( \frac{q - 1}{q}, \frac{1}{q (q - 1)},
     \cdots, \frac{1}{q (q - 1)} \right) . \]
  We prove now that the probability of having a giant component in $G (n, p_n,
  q)$ has the same limit as the probability
  \begin{equation}
    \mu_{\beta_n, h_n, n} \left( \frac{\tmmathbf{N}}{n} \notin B
    (\tmmathbf{x}^s, \varepsilon) \right) \label{eq-probaGC}
  \end{equation}
  for small enough $\varepsilon > 0$, for $\beta_n$ satisfying $p_n = 1 - \exp
  (- \beta_n / n)$ and $h_n = 0$.
  
  Indeed, let us fix a realization of the spins. Then we open edges between
  spins of equal color with probability $p_n$, resulting in a collection of
  $q$ Erd\"os-R\'enyi random graphs $G (N_i, p_n, 1)$ for $i = 1, \ldots, q$.
  It is known that a giant cluster appears in such a graph when $\lim_n N_i
  p_n > 1$ (see for instance {\cite{Bo}}). Yet, in the symmetric state one has
  $\lim_n N_i p_n = \beta_c / q < 1$ as $q > 2$, hence no giant component
  appears. In the asymmetric state $\tmmathbf{x}^{a, i}$ on the opposite, one
  has $\lim p_n N_i = \beta_c (q - 1) / q > 1$ thus a giant component emerges
  with conditional probability going to $1$.
  
  Finally, the quantity (\ref{eq-probaGC}) is computed using Theorem
  \ref{thm-probacoex} after we remark that
  \[ \beta_n = \beta_c + \frac{1}{n}  \left( \gamma + \frac{\beta_c^2}{2}
     \right) + o \left( \frac{1}{n} \right) . \]
\end{proof}

Let us conclude on the computation of the partition function for the
random-cluster model:

\begin{proof}
  (Proposition \ref{prop-Zrc}). We begin with a computation that permit to
  relate the partition function of the Curie-Weiss Potts model to that of the
  random-cluster model, defined at (\ref{eq-def-ZRC}). Now we say that an edge
  configuration $\omega \in \{0, 1\}^{E (K_n)}$ and a spin configurations
  $\sigma \in \{1, \ldots, q\}^n$ are {\tmem{compatible}} when $\omega_e = 1
  \Rightarrow \sigma_i = \sigma_j$, for all $e =\{i, j\} \in E (K_n)$, which
  we denote as $\omega \prec \sigma$. The factor $q^{C (\omega)}$ can be
  understood as the number of spin configurations $\sigma$ that are compatible
  with $\omega$. Hence:
  \begin{eqnarray}
    Z^{\tmop{RC}}_{p, q, n} & = & \sum_{\omega \in \{0, 1\}^{E (K_n)}}
    \sum_{\sigma \in \{1, \ldots, q\}^n : \omega \prec \sigma} \prod_{e \in
    K_n} p^{\omega_e} (1 - p)^{1 - \omega_e}  \nonumber\\
    & = & \sum_{\sigma \in \{1, \ldots, q\}^n} \sum_{\omega \in \{0, 1\}^{E
    (K_n)} : \omega \prec \sigma} \prod_{e \in K_n} p^{\omega_e} (1 - p)^{1 -
    \omega_e} \nonumber\\
    & = & \sum_{\sigma \in \{1, \ldots, q\}^n} \exp \left( - \frac{\beta}{n} 
    \sum_{1 \leqslant i < j \leqslant n} (1 - \delta_{\sigma_i, \sigma_j})
    \right) \nonumber\\
    & = & \sum_{\sigma \in \{1, \ldots, q\}^n} \exp \left( - \frac{\beta}{2}
    (n - 1) + \frac{\beta}{n}  \sum_{1 \leqslant i < j \leqslant n}
    \delta_{\sigma_i, \sigma_j} \right) \nonumber\\
    & = & Z_{\beta, 0, n} \exp \left( - \frac{\beta}{2} (n - 1) \right) . 
    \label{eq-ZRC-ZPotts}
  \end{eqnarray}
  for $\beta$ such that $p = 1 - \exp (- \beta / n)$. Remains to determine the
  asymptotics of $Z_{\beta_n, 0, n}$ for $\beta < \beta_c$. Thanks to the
  assumption $\beta < \beta_c$ the minimizer of the free energy is unique and
  symmetric:
  \[ \tmmathbf{x}^s = \left( \frac{1}{q}, \cdots, \frac{1}{q} \right) . \]
  This implies $\tmmathbf{d}_n = 0$ (see Remark \ref{rmk-d}), thus
  $\tmmathbf{N}= n\tmmathbf{x}^s + n^{1 / 2} \tmmathbf{W}$. Equation
  (\ref{equiv-Zmug}), in the limit $R \rightarrow \infty$, gives
  \begin{eqnarray}
    Z_{\beta_n, 0, n}  & = & (1 + o_n (1)) A_{\beta} (\tmmathbf{x}^s) e^{-
    nf_{\beta_n, 0} \left( \tmmathbf{x}^s \right)}  \int_{\mathcal{H}} e^{-
    \frac{1}{2} \tmmathbf{Q}_{\tmmathbf{x}^s, \beta} (\tmmathbf{w})^{}} \mathd
    w_2 \ldots \mathd w_q  \label{eq-Zmu}
  \end{eqnarray}
  as $\int_P \mathd w_2 \ldots \mathd w_q = 1$ for $\tmmathbf{w}= (- \sum_{i =
  2}^q w_i, w_2, \ldots, w_q)$. According to (\ref{eq-Aint}) one has
  \begin{eqnarray*}
    A_{\beta} (\tmmathbf{x}^s) \int_{\mathcal{H}} e^{- \frac{1}{2}
    \tmmathbf{Q}_{\tmmathbf{x}^s, \beta} (\tmmathbf{w})^{}} \mathd w_2 \ldots
    \mathd w_q & = & e^{- \frac{\beta}{2}}  \left( 1 - \frac{\beta}{q}
    \right)^{- \frac{q - 1}{2}} .
  \end{eqnarray*}
  On the other hand, the free energy is easily computed:
  \begin{eqnarray*}
    f_{\beta_n, 0} (\tmmathbf{x}^s) & = & \log \frac{1}{q} - \frac{\beta_n}{2
    q}
  \end{eqnarray*}
  leading to
  \begin{eqnarray*}
    Z_{\beta_n, 0, n} & = & (1 + o_n (1)) \left( 1 - \frac{\beta}{q}
    \right)^{- \frac{q - 1}{2}} q^n \exp \left( \frac{n \beta_n}{2 q} -
    \frac{\beta}{2} \right) .
  \end{eqnarray*}
  Then (\ref{eq-ZRC-ZPotts}) for $p_n = 1 - \exp (- \beta_n / n)$ gives:
  \begin{eqnarray*}
    Z^{\tmop{RC}}_{p_n, q, n} & = & Z_{\beta_n, 0, n} \exp \left( -
    \frac{\beta_n}{2} (n - 1) \right)\\
    & = & (1 + o_n (1)) \left( 1 - \frac{\beta}{q} \right)^{- \frac{q -
    1}{2}} q^n \exp \left( - \frac{n \beta_n}{2} \left( \frac{q - 1}{q}
    \right) \right)
  \end{eqnarray*}
  and the proof is over as $p_n = \beta / n + \gamma / n^2 + o (1 / n^2)$
  implies $\beta_n = \beta + (\gamma + \beta^2 / 2) / n + o (1 / n)$.
\end{proof}

\end{document}